\documentclass{llncs}
\usepackage[english]{babel}
\usepackage[pdftex]{hyperref}
\usepackage{algorithm}
\usepackage{algorithmic}
\usepackage[table]{xcolor}
\definecolor{red}{rgb}{0,0,0} 

\usepackage{graphicx}
\usepackage{amssymb, amsmath}
\usepackage{comment}
\usepackage{color}
\usepackage{subfigure}

\begin{document}

\title{The $\kappa_r$-version of the WRT-invariants, monochromatic 3-connected blinks 
and evidence for a conjecture on their induced 3-manifolds}

\date{\today}

\author{S\'ostenes Lins\inst{1}\thanks{Supported by CIn/UFPE and CNPq, Brazil.} \and Craig Hodgson\inst{3}\thanks{Supported by the Australian Research Council} \and Lauro Lins\inst{2} \and Cristiana Huiban\inst{1}\thanks{Supported by FACEPE and CNPq, Brazil.}}
\authorrunning{S.\ Lins, C.\ Hodgson, L.\ Lins, C.\ Huiban} 
\institute{Universidade Federal de Pernambuco, Brazil\\
\email{\{sostenes,cmngh\}@cin.ufpe.br}\\ 
\and
AT\&T Labs Research \\
\email{llins@research.att.com}\\
\and
University of Melbourne, Australia\\
\email{craigdh@unimelb.edu.au}
}
\maketitle

\begin{abstract}
A {\em blink} is a plane graph with a bipartition (black, gray) of its edges.
Subtle classes of blinks are in 1-1 correspondence with closed, oriented and connected 3-manifolds up to orientation preserving homeomorphisms~\cite{lins2013B}. Switching black and gray in a blink $B$, giving $-B$, reverses the manifold orientation.
The dual of the blink $B$ in the sphere $\mathbb{S}^2$ is denoted by $B^ \star$. Blinks $B$ and $-B^\star$ induce the same 3-manifold. 
The paper reinforces the Conjecture that if $B' \notin \{B,-B^\star\}$, then the monochromatic 3-connected (mono3c) blinks $B$ and $B'$ induce distinct 3-manifolds. Using homology of covers and length spectra, we conclude the topological classification of 708 mono3c blinks that were organized in equivalence classes by WRT-invariants in~\cite{lins2007blink}. We also present a reformulation of the combinatorial algorithm to obtain the WRT-invariants of~\cite{lins1995gca} using only the blink.
\end{abstract}

\section{Introduction}

It is well known that a 3-connected planar graph has a unique pair of embeddings in the 2-sphere up to orientation preserving homeomorphism~\cite{whitney1933}, see also~\cite{welsh2010matroid}. {\color{red}In other words, every planar 3-connected graph has only two planar maps, where their faces differ only in the orientation.} The present paper suggests a manifestation of this fact in closed oriented connected 3-manifolds induced by mono3c blinks, restricting the 1-1 correspondence of~\cite{lins2013B}~\href{http://arxiv.org/abs/1305.4540}{\underline{arxiv}} to these blinks. A blink is a finite plane graph with an arbitrary bipartition (black, gray) of its edges. 
Blinks provide a universal language for 3-manifolds~\cite{lins2013C}~\href{http://arxiv.org/abs/1305.5590}{\underline{arxiv}}, \cite{lins2013B}~\href{http://arxiv.org/abs/1305.4540}{\underline{arxiv}}.

\subsection{Blinks and $3$-manifolds}

\begin{figure}
\begin{center}
\includegraphics[width=9.0cm]{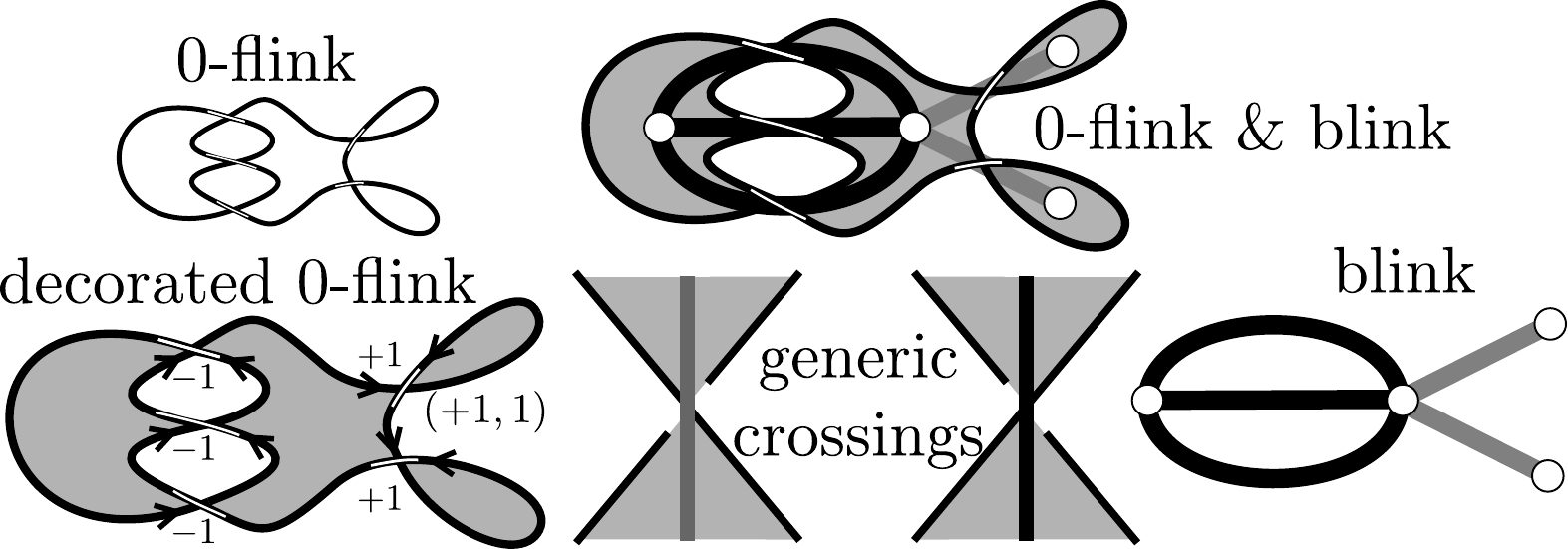} 
 \caption{\sf From 0-flink (i.e. a blackboard framed link, that is a projection of the framed link without twists) to blink and back: the projection of any 0-flink can be 2-face colored
into white and gray with the infinite face being white so that each sub-curve between two crossings has their incident faces receiving distinct colors. 
The above figure shows how to transform a 0-flink projection into a blink (with thicker edges than the curves representing the link projection). The vertices of the blink are distinguished 
points represented by white disks in the interior of the gray faces. Each crossing of the link projection becomes an edge in the corresponding blink. An edge of the blink is gray if the upper strand that crosses it is from northwest to southeast, it is black if the upper strand that crosses it is from northeast to southwest. The inverse procedure is clearly defined. In fact, the 0-flink is the so called {\em medial map} of the blink. Thus we have a 1-1 correspondence between 0-flinks and blinks. The flinks (non-negative fraction decorated links),  defined and treated in~\cite{lins2013B}~\href{http://arxiv.org/abs/1305.4540}{\underline{arxiv}}, are a generalization of blackboard framed links~\cite{kauffman1991knots,kauffman1994tlr}.}
\label{fig:flinkblink}
\end{center}
\end{figure}

\begin{figure}
\begin{center}
\includegraphics[width=8.0cm]{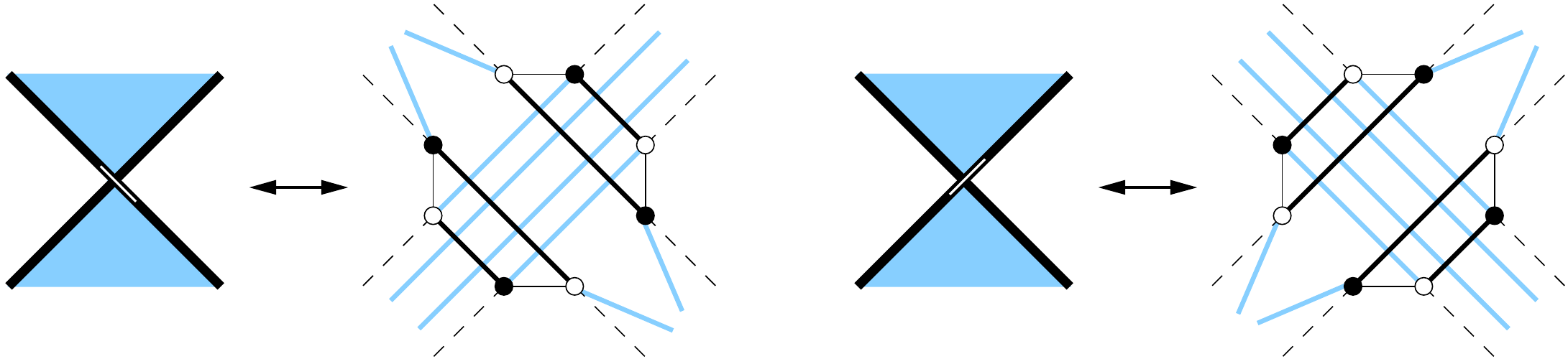} 
\caption{\sf Each crossing of the link is replaced with a subpart of the $3$-gem.}
\label{fig:gemblink}
\end{center}
\end{figure}

{\color{red}To understand how to obtain a $3$-manifold from a blink we need the following definitions. An $n$-residue is a connected component of a edge-colored graph induced by all the edges of $n$ chosen colors. A \emph{$3$-gem} is a finite $4$-regular graph with a $4$-proper edge coloring in which $v + t = b$, where $v$ is the number of vertices, $b$ is the number of $2$-residues and $t$ is the number of $3$-residues. It follows from the \emph{Triangulation Theorem for 3-manifolds} of Moise~\cite{moise1952affine} that every closed compact (as oriented and connected) $3$-manifold can be induced by a $3$-gem. 

A $k$-simplex is a k-dimensional polytope which is the convex hull of its $k + 1$ vertices (that is $C=\{\phi_0 u_o+ ...+ \phi_k u_k|\phi_i\geq0, 0\leq i \leq k, \sum_{i=0}^{k}\phi_i=1\}$ with $u_0, ..., u_k \in \mathbb{R}^k$ affinely independent). A simplicial complex $\mathcal{K}$ is a set of simplices where any face of a simplex from $\mathcal{K}$ is also in $\mathcal{K}$ and, the intersection of any two simplices $\sigma_1, \sigma_2 \in \mathcal{K}$ is a face of both $\sigma_1$ and $\sigma_2$. A triangulation of a topological space $X$ (as a $3$-manifold) is a simplicial complex $\mathcal{K}$ with a homeomorphism (i.e. a bicontinuous function) $h:\mathcal{K} \to X$. 

To obtain a triangulation from a blink, we make the transformations $blink \rightarrow link \rightarrow$ $3$-gem $\rightarrow triangulation$. Fig.~\ref{fig:flinkblink} shows how to obtain a blink from a \emph{link}, and vice versa. Fig.~\ref{fig:gemblink} describes the transformations between \emph{link} and \emph{3-gem}~\cite{lins1995gca,lins2007blink}.}

{\color{red}To obtain a triangulation from a $3$-gem, create a $d$-point for each vertex (d=0), edge (d=1), 2-residue (d=2) and 3-residue (d=3) of the $3$-gem. To obtain the space coordinates of these points, we define a code as follows. Assign the code $k^c$ for each $0$-point, where $k$ is an integer number and $c=\{\}$. For each $1$-point, $k$ is the smallest $k$ among the neighbors and $c$ is the color of its corresponding edge in the $3$-gem. For each $2$-point, $k$ is also the smallest $k$ assigned to the neighbors and $c$ is the union of sets $c$ of all neighbors. So, an example of tetrahedra formed is $\{k, k^1, k^{12}, k^{123}\}$.

In order to define a valid simplicial complex, we need to obtain an arrangement of simplices such that the intersections are only a point, a side or a face of the defined simplices. To guarantee that two $1$-coordinate point are disjoint we need $\mathbb{R}^1$, for 2 line segments (4 points) we need $\mathbb{R}^3$, for our case that is for two tetrahedras (eight points) we need $\mathbb{R}^7$. The next step is then to embed the points into $\mathbb{R}^7$. Each code $k^c$ is mapped to an integer number $i$ and its coordinate in $\mathbb{R}^7$ will be $\{i^1,i^2,\dots, i^7\}$.

To obtain an $3$-gem back from a given triangulation, apply the barycentric subdivision for each tetrahedra. So, each tetrahedra is subdivided in $24$ sub-tetrahedras, $6$ per face. Assign color $k$ to each $k$-simplex and color the faces of each sub-tetrahedra with the color of its opposite vertex. Each sub-tetrahedra becomes a vertex of the $3$-gem. Add an edge between two vertices coming from two sub-tetrahedras sharing $3$ vertices (defining a face $F$) and assign to the edge the color of $F$.

From $3$-gem to blink, see~\cite{linsMachado2013} for an $O(n^2)$ algorithm (where $n$ is the number of tetrahedras) to obtain a framed link from a special type of triangulation. The general transformation is an open subject.
}

\subsection{$HGnQI$-classes and the topological classification of $T_{16}$}
\begin{figure}
\begin{center}
\includegraphics[width=12.0cm]{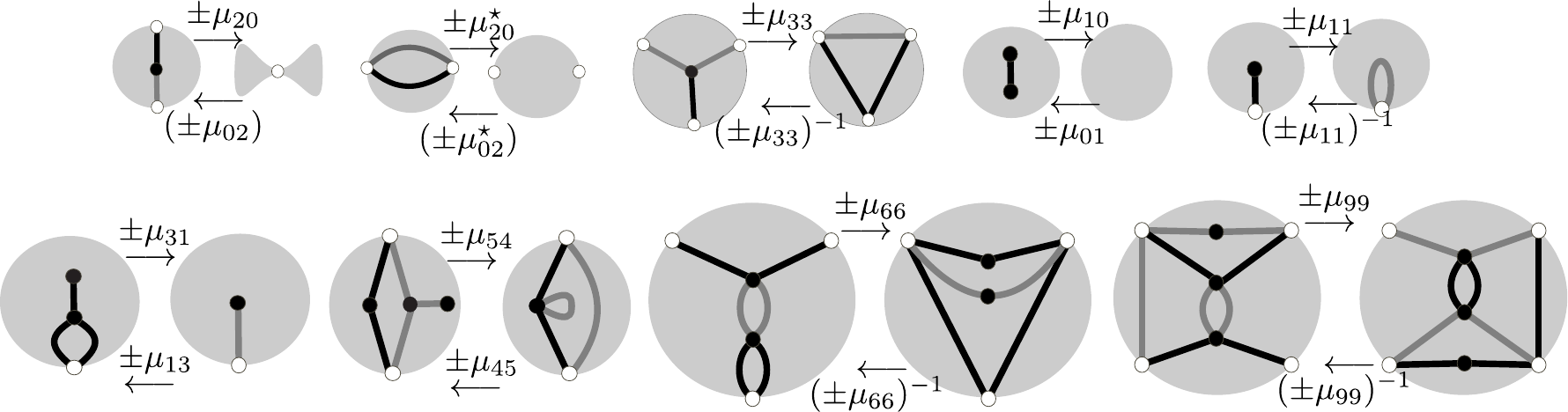} 
\caption{\sf The set of blinks inducing the $3$-manifolds from the same homeomorphism class is closed under these blink transformations called \emph{coin calculus}. The complete set of transformations has 36 moves (not necessarily distinct) that are obtained from the moves in the figure by dualizing and interchanging black and gray edges.}
\label{fig:reducedCoinCalculus}
\end{center}
\end{figure}

An $HGnQI$-class of blinks is the set of blinks inducing a 3-manifold with the same Homology Group and the same $\kappa_r$ (Quantum) Invariants for $r=3,\ldots,n$. The latter corresponds to the Kauffman-Lins version~\cite{kauffman1994tlr} of the WRT-invariants (Witten-Reshetikhin-Turaev invariants). {\color{red} If $B'$ is another blink obtained from $B$ by moves in the coin calculus (see Fig.~\ref{fig:reducedCoinCalculus}), they induce the same $3$-manifold  up to orientable homeomorphisms. Otherwise, $\kappa_r(B) \neq \kappa_r(B')$.}

Our idea in Section~\ref{section:kappaerre} is to make more available for the Combinatorics community these strong and mysterious quantum invariants: an infinite sequence of complex numbers which can be deterministically assigned to blinks up to \emph{coin calculus}~\cite{lins2013B}~\href{http://arxiv.org/abs/1305.4540}{\underline{arxiv}} and to 3-manifolds up to orientable homeomorphisms. The invariance issue will remain untouched here; it is treated in detail in~\cite{kauffman1994tlr}. 

{\color{red} Let $T_{16}$ (all mono3c blinks up to 16 edges) be the set of 708 mono3c blinks distributed in 381 $HG8QI$ classes by WRT-invariants, as described  in~\cite{lins2007blink}~\href{http://arxiv.org/abs/math/0702057}{\underline{arxiv}}.} 
The set was explicitly generated, displayed and classified up to oriented homeomorphisms (for the 3-manifolds induced by the blinks) by, $\kappa_r$, $r=3,\ldots 8$, leaving exactly 11 doubts. These doubts correspond to eleven $HG8QI$-classes with more than one pair of blinks $(B, -B^*)$. We complete the classification by showing that they give non-homeomorphic pairs in Section~\ref{Sec:HG8QI_classification}.

The topological classification of $T_{16}$ reinforces the Conjecture in~\cite{lins2013B} that if $B' \notin \{B,-B^\star\}$, then the mono3c blinks $B$ and $B'$ induce distinct 3-manifolds. These experimental results are  evidence that the oriented homeomorphism problem for those manifolds is solvable in $O(n^2\log n)$-time, where $n$ is the number of edges in the blinks inducing them. This is the complexity of the isomorphism problem for blinks~\cite{lins1980sequence}.

In Section~\ref{Sec:HG8QI_classification}, we were able to distinguish all the remaining $HG8QI$-classes by  the lengths of closed geodesics in their hyperbolic structures, using the method presented in~\cite{hodgson1994symmetries}. The running time was substantially better in some cases compared with the times obtained with the homology of covers technique also used in this section. We have used a combination of the software SnapPy~\cite{snappy}, GAP~\cite{gap} and Sage~\cite{sage}
to make the classification. In Section~\ref{section:kappaerre}, for completeness, we present the full definition and a reformulation of the combinatorial algorithm to obtain the WRT-invariants of~\cite{lins1995gca} using only the blink.

\section{$\kappa_r$: the Kauffman-Lins version of the WRT-invariants}
\label{section:kappaerre}

The invariant $\kappa_r$ was obtained and justified in the Kauffman-Lins monograph \cite{kauffman1994tlr}. Given an integer number $r \ge 3$, a complex number $\kappa_r(B)$ is associated to each blink $B$. If $B'$ is another blink obtained from $B$ by moves in the coin calculus, then $\kappa_r(B)=\kappa_r(B')$. Therefore, $\kappa_r$ is an {\em invariant of closed, orientable, connected 3-manifolds}. 

This section is a reformulation of Chapter 7 of~\cite{lins1995gca}. The invariance of $\kappa_r$ relies on the properties of Kauffman's bracket and on the algebraic properties of the Temperley-Lieb algebra. Subsequently it was realized that this invariant is one of the manifestations of the Witten-Reshetikhin-Turaev invariants.
The complete theory is developed from scratch in~\cite{kauffman1994tlr}.

Originally these invariants were found by Witten in the late 1980's using a physical formalism that was not mathematically satisfactory. Witten's result broke the prejudice that good invariants of 3-manifolds did not exist. Shortly after, some eastern European researchers such as Turaev, Viro, Reshetikhin, Kirillov and others, produced full mathematical proofs of the invariance of Witten's results using quantum groups~\cite{witten1989quantum,reshetikhin1991invariants,turaev1992state,kirillov1993orbit}.
However, quantum groups is a rather complicated subject, so Kauffman-Lins makes possible to combinatorialize the whole situation. They used some ideas of Lickorish (invariance under the second of Kirby's moves~\cite{lickorish1991three}) via the Temperley-Lieb algebra and cubic graphs embedded into 3-space, providing the $\kappa_r$-invariant. It demands much less machinery to be understood--- see also page 144 of Kauffman-Lins monograph, and Turaev's \emph{shadow world} in~\cite{turaev1994quantum}.

In Subsection~\ref{Subsec:Algebraic}, we describe the algebraic ingredients to compute the factors of the blink (in Subsection~\ref{Subsec:Factors}) that are used to compute the invariant. It is our intention here to provide a complete recipe for obtaining the $\kappa_r$-invariant for blinks in a way to be simply understood by a combinatorially inclined reader.

\subsection{Algebraic ingredients}
\label{Subsec:Algebraic}

In this subsection, we present all the algebraic ingredients we need to define the function $\kappa_r$ (for a fixed integer $r \ge 3$) on a blink $B$. Let $A=e^{\frac{\pi i}{2r}}$ be the ``first'' primitive $4r$-th root of unity and ${\cal I} = \{0,1,...,r-2\}$. For $n$ in ${\cal I}$, let $$\Delta_n = (-1)^n\frac{\textstyle A^{2n+2}-A^{-2n-2}} {\textstyle A^2-A^{-2}}, \hspace{10mm} [n] = \frac{\textstyle A^{2n}-A^{-2n}}{\textstyle A^2-A^{-2}}=(-1)^{n-1}\Delta_{n-1}.$$ Letting $q=A^2$, for reasons inherited from the physics we call $[n]$ the {\em $q$-deformed quantum integer}\index{$q$-deformed!quantum integer} and $[n]! = \prod_{0 \le m \le n} [m]$ the {\em $q$-deformed quantum factorial}\index{$q$-deformed!quantum factorial}. Note that even though $A$ is complex, $\Delta_n$ and $[n]$ assume only real values. 

Three numbers $a,b,c \in {\cal I}$ form an {\em admissible triple}\index{admissible triple} if $a+b+c \le 2r-4$ and $a+b-c, b+c-a, c+a-b$ are non-negative and even. Let $\theta: \mathcal{I}^3 \mapsto \mathbb{R}$, defined on the admissible triples by means of
$$
\theta(a,b,c)=\frac{\textstyle (-1)^{m+n+p}[m+n+p+1]! [n]! [m]! [p]!}
{\textstyle [m+n]! [n+p]! [p+m]!}\ , \mbox{where}
$$ 
$m= (a+b-c)/2, n=(b+c-a)/2, p=(c+a-b)/2.$
We define $\theta(a,b,c)=0$ if $(a,b,c)$ fails to be admissible.

Let $\lambda: \mathcal{I}^3 \mapsto \mathbb{C}$, be defined on the admissible triples by
$$
\lambda(a,b,c)=\lambda^{ab}_c=
(-1)^{(a+b-c)/2}A^{[a(a+2)+b(b+2)-c(c+2)]/2}
$$

We let $\lambda^{ab}_c = 0$ if $(a,b,c)$ fails to be admissible. Finally, define $Tet: \mathcal{I}^6 \mapsto \mathbb{R}$, as follows. If $(\alpha,\beta,\phi)$, $(\alpha,\delta,\epsilon)$, $(\gamma,\delta,\phi)$, $(\beta,\gamma,\epsilon)$ are admissible, define

$$Tet(\alpha,\beta,\gamma,\delta,\epsilon,\phi)=
Tet\left[
\begin{matrix}
\alpha&\beta&\epsilon \cr 
\gamma&\delta&\phi\cr
\end{matrix}
\right]
=\frac{\mathcal
{I}nt!}{\mathcal{\epsilon}xt!}\sum_{m\leq s\leq
M}\frac{(-1)^s[s+1]!}{\prod_i[s-a_i]!\prod_j[b_j-s]!}$$ 
where,
$$\begin{array}{rcl}
{\cal I}nt!&=&\prod_{i,j}[b_j-a_i]!\\ [0.2cm]
{\cal \epsilon}xt!&=&[\alpha]![\beta]![\gamma]![\delta]![\epsilon]![\phi]!\\ [0.2cm]
a_1&=&\frac{1}{2}(\alpha+\delta+\epsilon)\qquad b_1=\frac{1}{2}(\beta+\delta+\epsilon+\phi)\\ [0.2cm]
a_2&=&\frac{1}{2}(\beta+\gamma+\epsilon)\qquad b_2=\frac{1}{2}(\alpha+\gamma+\epsilon+\phi)\\ [0.2cm]
a_3&=&\frac{1}{2}(\alpha+\beta+\phi)\qquad b_3=\frac{1}{2}(\alpha+\beta+\gamma+\delta)\\ [0.2cm]
a_4&=&\frac{1}{2}(\gamma+\delta+\phi)\qquad m\!=\!\max\{a_i\}\quad M\!=\!\min\{b_j\}\\
\end{array}$$
If one of the four triples above fails to be admissible, define the value of $Tet$ as null. Note that the bipartition of the edge set of the blink is disregarded, except for the $\lambda^{ab}_c$'s. They are also the only terms which are indeed complex: $Tet(\alpha,\beta,\gamma,\delta,\epsilon,\phi)$, $\theta(a,b,c)$ and $\Delta_n$ are reals.

See Algorithm~\ref{alg:Algebraic} for a summary of the subsection.

\begin{algorithm}
  \caption{Computation of the algebraic ingredients}
  \label{alg:Algebraic}
  \begin{algorithmic}[1]
    \STATE{$A=e^{\frac{\pi i}{2r}}$, $q=A^2$, ${\cal I} = \{0,1,...,r-2\}$} 
    \FOR{All $n$ in ${\cal I}$}
	\STATE{$\Delta_n = (-1)^n\frac{\textstyle A^{2n+2}-A^{-2n-2}}{\textstyle A^2-A^{-2}}$}
	\STATE{$[n]=(-1)^{n-1}\Delta_{n-1}$}
	\STATE{$[n]! = \prod_{0 \le m \le n} [m]$}
    \ENDFOR
    \FOR{All $(a,b,c) \in {\cal I}^3$}
    \IF{$(a,b,c)$ is \em admissible}
	\STATE{$m= (a+b-c)/2$, $n=(b+c-a)/2$, $p=(c+a-b)/2$}
        \STATE{$\theta(a,b,c)=\frac{\textstyle (-1)^{m+n+p}[m+n+p+1]! [n]! [m]! [p]!}{\textstyle [m+n]! [n+p]! [p+m]!}$}
        \STATE{$\lambda^{ab}_c=(-1)^{(a+b-c)/2}A^{[a(a+2)+b(b+2)-c(c+2)]/2}$}       
    \ELSE
      \STATE{$\theta(a,b,c)=0$, $\lambda^{ab}_c=0$}
    \ENDIF
    \ENDFOR
    \FOR{All $(\alpha,\beta,\gamma,\delta,\epsilon,\phi) \in \mathcal{I}^6$}
      \IF{$(\alpha,\beta,\phi)$,$(\alpha,\delta,\epsilon)$,$(\gamma,\delta,\phi)$,$(\beta,\gamma,\epsilon)$ are \textit{all} admissible}
        \STATE{${\cal I}nt!=\prod_{i,j}[b_j-a_i]!$}
	\STATE{${\cal \epsilon}xt!=[\alpha]![\beta]![\gamma]![\delta]![\epsilon]![\phi]!$}
	\STATE{$a_1=\frac{1}{2}(\alpha+\delta+\epsilon)$, $a_2=\frac{1}{2}(\beta+\gamma+\epsilon)$}
	\STATE{$a_3=\frac{1}{2}(\alpha+\beta+\phi)$, $a_4=\frac{1}{2}(\gamma+\delta+\phi)$} 
	\STATE{$b_1=\frac{1}{2}(\beta+\delta+\epsilon+\phi)$, $b_2=\frac{1}{2}(\alpha+\gamma+\epsilon+\phi)$, $b_3=\frac{1}{2}(\alpha+\beta+\gamma+\delta)$}
        \STATE{$m\!=\!\max\{a_i\}$,$M\!=\!\min\{b_j\}$}
	\STATE{$Tet\left[
	  \begin{matrix}
	    \alpha&\beta&\epsilon \cr 
	    \gamma&\delta&\phi\cr
	  \end{matrix}
	\right]
	=\frac{\mathcal
	{I}nt!}{\mathcal{\epsilon}xt!}\sum_{m\leq s\leq
	M}\frac{(-1)^s[s+1]!}{\prod_i[s-a_i]!\prod_j[b_j-s]!}$}
      \ELSE \STATE{$Tet\left[
	  \begin{matrix}
	    \alpha&\beta&\epsilon \cr 
	    \gamma&\delta&\phi\cr
	  \end{matrix}
	\right]=0$}
      \ENDIF
   \ENDFOR
  \end{algorithmic}
\end{algorithm}

\subsection{Blink factors and the $\kappa_r$-invariant definition}
\label{Subsec:Factors}

\begin{figure}
\begin{center}
\includegraphics[width=4cm]{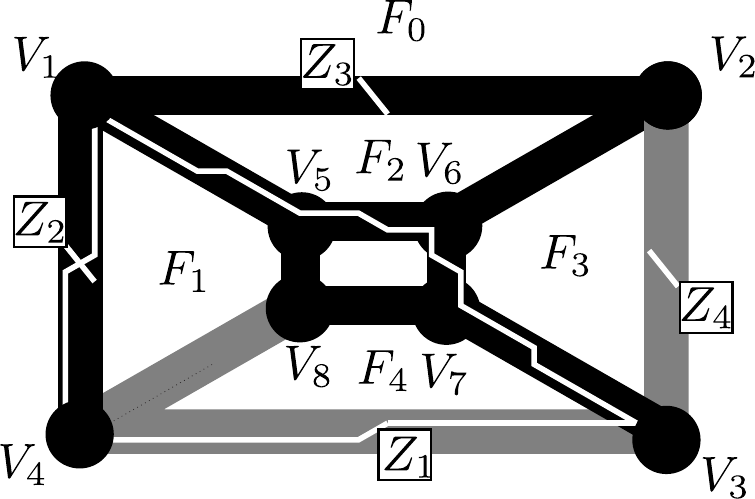} 
\caption{\sf A \emph{zigzag} in blink $B$ (as in any planar graph) is a closed path that alternates taking the leftmost turn and the rightmost turn at each vertex; zigzags in general surfaces (orientable and non-orientable) are studied in \cite{lins1982graph}. The bipartition of the edge colors of the blink plays no role in the definition of the zigzags; a zigzag traverses each edge at most twice, so that the total number of edges traversed by the set of zigzags is exactly $2n$, where $n$ is the number of edges of the blink. An \emph{angle} is a pair of consecutive edges in the boundary of a face $F$. Note that the edges of an angle are consecutive in the edge-star of a vertex $V$ and in a zigzag path $Z$. So an angle defines a triple $(V,F,Z)$. In the example, 
the zigzag $Z_1$ is shown and the domain of $\sigma$ is $V_B \cup F_B \cup Z_B=\{V_1,V_2,V_3,V_4,V_5,V_6,V_7,V_8,F_0,F_1,F_2,F_3,F_4,F_5,Z_1,Z_2,Z_3,Z_4\}$.
}
\label{fig:zigzagcube}
\end{center}
\end{figure}

{\em Factors} (complex numbers) are associated to the following {\em sites} of $B$: black edge, gray edge, angle, vertex, face and zigzag (see definitions in Fig.~\ref{fig:zigzagcube}). To compute the factors, we use the algebraic ingredients defined before: $Tet(\alpha,\beta,\gamma,\delta,\epsilon,\phi)$, $\lambda^{ab}_c$, $\overline{\lambda^{ab}_{c}}$ (multiplicative inverse of $\lambda^{ab}_c$),  $\theta(a,b,c)$, $\Delta_n$.

Let $\sum = \{\sigma:\{\mathcal{V_B} \cup \mathcal{F_B} \cup \mathcal{Z_B}\} \mapsto \mathcal{I}\}$ be the set of all possible \emph{blink states}, where $\{\mathcal{V_B} \cup \mathcal{F_B} \cup \mathcal{Z_B}\}$ is the union of the vertices, faces and zigzags of blink $B$. Each state $\sigma$ is a specific mapping that associates to each vertex, face and zigzag of $B$ a number in $\mathcal{I}$. Taking the trio $(V,F,Z)$ representing an angle of $B$, we generate the triple $(\sigma(V),\sigma(F),\sigma(Z))=(v,f,z)$. We say that a blink state is \emph{admissible} if each trio of $B$ has an admissible triple. If a state is not admissible, its value is declared to be 0. 

\begin{figure}
\begin{center}
\includegraphics[angle=-90,width=12.0cm]{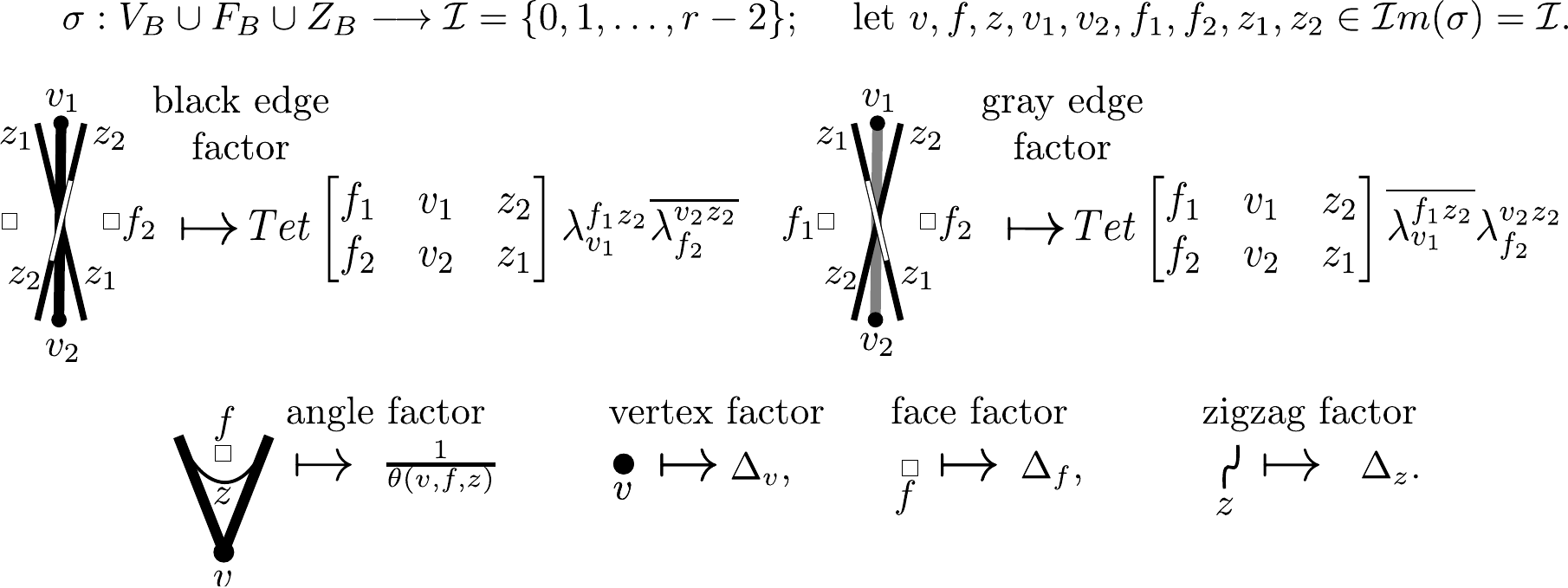} 
\caption{\sf $v,f,z,v_1,v_2,f_1,f_2,z_1,z_2$ are the images under the specific state $\sigma$ of, respectively, $V,F,Z,V_1,V_2,F_1,F_2,Z_1,Z_2$. There are two zigzags crossing each edge of the blink. So we represent each edge by sextuplets $(v_1, f_1, z_1, v_2, f_2, z_2)$. The functions used on the factors computation, $Tet(\alpha,\beta,\gamma,\delta,\epsilon,\phi)$, $\lambda^{ab}_c$, $\overline{\lambda^{ab}_c}$ (multiplicative inverse of $\lambda^{ab}_c$),  $\theta(a,b,c)$, $\Delta_n$ are detailed in Subsection~\ref{Subsec:Algebraic}. There are symmetries in the situation which imply that the edge factors are invariant under a half turn rotation. }
\label{fig:edgeanglefactors}
\end{center}
\end{figure}

The formulae for the factors computation of each site of $B$ are depicted in Fig.~\ref{fig:edgeanglefactors}. The product of these factors defines \emph{the value} of state $\sigma$, let us say $[\sigma]$. Let $\kappa'_r(B)$ be the sum of the values of all states in $\sum$.

The value $\kappa'_r(B)$ turns out to be an invariant under Kirby's second move, the {\em band move} or the {\em handle slide move}~\cite{kirby1978calculus}, see also page 144 of~\cite{kauffman1994tlr} or~\cite{lins2013B}~\href{http://arxiv.org/abs/1305.4540}{\underline{arxiv}}. Having this invariance, all we need to have a 3-manifold invariant is to define

$\kappa_r(B) = \kappa_r'(B) \left(\sin(\frac{\pi}{r}) \sqrt{\frac{2}{r}}\right)^{|F|+1}
\left(
(-1)^{r-2} e^{\frac{3 i \pi(r-2)}{4r}}\right)^{-n(F)}.$

In this formula, $F$ represents the 0-flink associated to $B$ (see $F\leftrightarrows B$ transformations in Fig.~\ref{fig:flinkblink}), $|F|$ is the number of components of $F$, $n(F)$ is the number of positive eigenvalues minus the number of negative eigenvalues of the {\em linking matrix} of $F$ (see page 254 of~\cite{lins1995gca}). It is the symmetric matrix whose entries are the \emph{linking numbers} for the pairs of components of $F$. That is, for components $i$ and $j$, the linking number is $\mathcal{L}_{i,j} = \sum_{x \in \chi_{i,j}} \frac{1}{2} sg(x),$ where $\chi_{i,j}$ denotes the set of crossings of component $i$ and $j$, and $sg(x)$ is the sign of the crossing. 

Note that we introduce a factor to go from $\kappa_r'(B)$ to $\kappa_r(B)$. The factor is needed to keep the invariance under the Kirby's move of type 1 (that is when we introduce an isolated gray or black edge in the blink). With these definitions it is possible to prove (see page 146 in~\cite{kauffman1994tlr}) that for all $r \ge 3$,

$\kappa_r(\circ)=1 \rm{\ and\ }
\kappa_r(\circ \hspace{-0.5mm} \text{---} \hspace{-0.5mm} \circ)=\sin(\frac{\pi}{r}) \sqrt{\frac{2}{r}},$

\noindent where $\circ$ denotes the isolated vertex blink, inducing $\mathbb{S}^2\times \mathbb{S}^1$ and $\circ$\hspace{-0.5mm}---\hspace{-0.5mm}$ \circ$ denotes the isolated single edge blink inducing $\mathbb{S}^3$. A face of a connected blink is a connected component of complement $\mathbb{R}^2 \backslash B$. There is a distinguished face (the infinite one) which is mapped to 0 by any admissible state $\sigma$. With this constraint and the introduced factor, $\kappa_r'$ presents a multiplicative property, that is $\kappa'_r(B \cup B') = \kappa_r'(B) \times \kappa_r'(B')$ for disjoint blinks $B$ and $B'$.

See Algorithm~\ref{alg:Kappar} for the complete set of instructions to compute the $\kappa_r$-invariant. 
\begin{algorithm}
  \caption{$\kappa_r$: the Kauffman-Lins version of the WRT-invariants using only blinks \footnotesize{(Recall $\sigma(V)=v,\sigma(F)=f,\sigma(Z)=Z$, and see $F\leftrightarrows B$ transformations in Fig.~\ref{fig:flinkblink})}}
  \label{alg:Kappar}
  \begin{algorithmic}[1]
   \STATE{$\kappa_r'(B)=0$}
   \FOR{All $\sigma \in \sum$ of blink $B$}
     \FOR{All black edge $(V_1,F_1,Z_1,V_2,F_2,Z_2)$ of blink $B$}
      \STATE{$\Phi_{1}=Tet\left[
	  \begin{matrix}
	    f_1&v_1&z_2 \cr 
	    f_2&v_2&z_1\cr
	   \end{matrix}
	\right]\lambda^{f_1z_2}_{v_1}\overline{\lambda^{v_2z_2}_{f_2}}$}
      \ENDFOR
      \FOR{All gray edge $(V_1,F_1,Z_1,V_2,F_2,Z_2)$ of blink $B$}
	\STATE{$\Phi_{2}=Tet\left[
	  \begin{matrix}
	    v_1&v_1&z_2 \cr 
	    f_2&v_2&z_1\cr
	  \end{matrix}
	\right]\overline{\lambda^{f_1z_2}_{v_1}}\lambda^{v_2z_2}_{f_2}$}
      \ENDFOR
      \FOR{All angle $(V,F,Z)$ of blink $B$}      
	\STATE{$\Phi_{3}=\frac{1}{\theta(v,f,z)}$}
      \ENDFOR
      \FOR{All vertex $V$ of blink $B$}      
	\STATE{$\Phi_{4}=\Delta_{\sigma(V)}=\Delta_v$}
      \ENDFOR
      \FOR{All face $F$ of blink $B$}
	\STATE{$\Phi_{5}=\Delta_{\sigma(F)}=\Delta_f$}
      \ENDFOR
      \FOR{All zigzag $Z$ of blink $B$}      
	\STATE{$\Phi_{6}=\Delta_{\sigma(Z)}=\Delta_z$}
      \ENDFOR      
      \STATE{$[\sigma]=\prod_{i} \Phi_i$}
      \FOR{All angles $(V,F,Z)$ of blink $B$}
	  \IF{$(v,f,z)$ is not admissible}	  
	    \STATE{$[\sigma]=0$, \textbf{break for}}	      
	  \ENDIF
      \ENDFOR
      \STATE{$\kappa_r'(B)=\kappa_r'(B)+[\sigma]$}
   \ENDFOR
   \STATE{$|F|=$ \#components of 0-flink $F$ associated to blink $B$}
   \STATE{$n(F)=$\#positive eigenvalues$ - $\#negative eigenvalues of the {\em linking matrix} of $F$} 
   \STATE{$\kappa_r(B) = \kappa_r'(B) \left(\sin(\frac{\pi}{r}) \sqrt{\frac{2}{r}}\right)^{|F|+1}\left((-1)^{r-2} e^{\frac{3 i \pi(r-2)}{4r}}\right)^{-n(F)}.$}
  \end{algorithmic}
\end{algorithm}

\section{Classifying the $HG8QI$-classes as non-homeomorphic pairs}
\label{Sec:HG8QI_classification}

An $HG8QI$-class of blinks is the set of blinks inducing a 3-manifold with the same homology and the same $\kappa_r$ invariants, for $r=3,\ldots,8$. Even though they induce the same 3-manifold, dual blinks are not filtered in \cite{lins2007blink} for the generation of the $HG8QI$-classes (this was on purpose, to control the generation algorithms). So, these dual blinks $(\pm B, \mp B^\star)$ appear together in the same class. In fact we only draw one of $\pm B$ (see Appendix), the other is obtained by exchanging black and gray edges. The effect of taking the negative blink in $\kappa_r(B)$ is to conjugate the complex number: for all blinks $B$, $\kappa_r(B) = \overline{\kappa_r(-B)}$ and $\kappa_r(B) = \kappa_r(-B^\star)$.


From~\cite{lins2007blink}~\href{http://arxiv.org/abs/math/0702057}{\underline{arxiv}}, we know that among all 381 $HG8QI$-classes of $T_{16}$, there are only eleven with more than one pair of blinks $(B, -B^*)$. These are the following:
(1) $14^t_{24}$,\ 
(2) $15^t_{16}$,\  
(3) $15^t_{19}$,\  
(4) $15^t_{22}$,\  
(5) $16^t_{42}$,\  
(6) $16^t_{56}$,\  
(7) $16^t_{140}$,\  
(8) $16^t_{141}$,\  
(9) $16^t_{142}$,\  
(10) $16^t_{149}$,\  
(11) $16^t_{233}$.
All the other pairs of blinks $(B,-B^\star)$ in $T_{16}$ are complete graphical invariants for their induced 3-manifolds. We have used a combination of the softwares SnapPy, GAP and Sage to prove that all pairs of blinks (that are not $(\pm B, \mp B^\star)$) in these $HG8QI$-classes induce distinct 3-manifolds. 

Dual and negative pairs of monochromatic blinks $(\pm B, \mp B^\star)$ are in 1-1 correspondence with the 0-flinks (blackboard framed links), see the 0-flink $\leftrightarrows$ blink transformations in Fig.~\ref{fig:flinkblink}. The appendix presents both the figures of the blink and its associated 0-flink. However, the 0-flinks (good for input into SnapPy) and its pair of surgery coefficients attached to each component of the 0-flink are redundant, as they are implied by the small blink displayed at the southwest of each 0-flink. All the blinks are monochromatic, formed by black edges only. This corresponds to having only alternating 0-flinks. 

The 3-manifolds are distinguished by the combined application of SnapPy, GAP and Sage using calculations of homology for finite sheeted coverings in Subsection~\ref{subsec:covers} and using the geodesic length technique from~\cite{hodgson1994symmetries} in Subsection~\ref{subsec:geodesic}. The results are summarised in Table~\ref{Table:results}, and details of the homology and length spectra computations are given in the appendix.

In all the experiments, the triangulations (files with extension $.tri$) are created from a \emph{link diagram} that we draw with the link editor of SnapPy (we can access the editor typing $M=Manifold()$ on SnapPy). Then select \textit{Tools - Send to SnapPy} in the editor to load the link complement as the variable $M$. Finally we do Dehn fillings in SnapPy using coefficients $(a,1)$ and $(b,1)$, $M.dehn\_fill([(a,1),(b,1)])$. The values $a$, $b$ are the integer framings on the components corresponding to the self-writhe of the component for the blackboard framing. 

Note that SnapPy can also do rational Dehn fillings: if $p$,$q$ are relatively prime integers then $(p,q)$ filling in SnapPy corresponds to filling on the slope $\frac{p}{q}$, i.e. $p*(homological\_meridian) + q*(homological\_longitude)$ bounds a disk in the added solid torus. After Dehn filling, we check that $M.solution\_type()$ gives the output `all tetrahedra positively oriented' (to be sure we have the correct hyperbolic structure), and save this triangulation using $M.save()$. (If other output is obtained, we retriangulate with $M.randomize()$ before saving).

\subsection{Calculations of homology for finite sheeted coverings}
\label{subsec:covers}
We first find all index $k$ subgroups of the fundamental group of the induced 3-manifold  for small values of $k$, then compute the homology of the corresponding covers of the induced 3-manifolds.  In fact, using hyperbolic volumes and homology of covers as computed by SnapPy/GAP/Sage we were able to distinguish all the $HG8QI$-classes. The running time of SnapPy/Sage for the computation of covers varies a lot. In Table~\ref{Table:results}, see the homology of covers distinguishing the two manifolds and their computation time  on a LG ultrabook with1.70GHz Intel processor i5-3317U. 
To compute the volume/covers, see example (with degree $k=5$) of SnapPy/Sage script below.


{\footnotesize
\begin{verbatim}
sage: import snappy
sage: A=snappy.Manifold('T71.tri')
sage: B=snappy.Manifold('T79.tri')
sage: A.volume()
sage: B.volume()
sage: coversA=A.covers(5,method='gap')
sage: coversB=B.covers(5,method='gap')
sage: sorted[X.homology() for X in coversA]
sage: sorted[X.homology() for X in coversB]
\end{verbatim}
}

The computation of covers for two cases, $16^t_{140}$ and $16^t_{149}$, took more than three days using SnapPy/Sage and was inconclusive. The case of $16^t_{149}$ is already distinguished by the volumes. 

For the $HG8QI$-class $16^t_{140}$, there are no $k$-covers for $k=2,\ldots 6$ and the search for $7$-covers took more than a week without ending using SnapPy/Sage. So we did not try $k$-covers for $k \ge 8$. However, Nathan Dunfield observed that they can be distinguished using covers from representations onto $PSL(2,7)$ as computed using GAP/Sage, see script below. (Here the subgroups have index 8, corresponding to the stabilizer of a point under the natural action of $PSL(2,7)$ on $P^1({\mathbb F}_7)$.)

{\footnotesize
\begin{verbatim}
sage: import snappy
sage: M1 = snappy.Manifold('T423.tri')
sage: M2 = snappy.Manifold('T444.tri')
sage: G1=gap(M1.fundamental_group())
sage: G2=gap(M2.fundamental_group())
sage: Q=PSL(2,7)
sage: homs1,homs2=G1.GQuotients(Q),G2.GQuotients(Q)
sage: [M1.cover(h).homology() for h in homs1]
[Z/2 + Z/2 + Z/4 + Z/4 + Z/641184, Z/2 + Z/2 + Z/4 + Z/4 + Z/641184]
sage: [M2.cover(h).homology() for h in homs2]
[Z/2 + Z/2 + Z/4 + Z/4 + Z/221232, Z/2 + Z/2 + Z/4 + Z/4 + Z/774480]
\end{verbatim}
}

\subsection{Geodesic length technique}
\label{subsec:geodesic}
All the pairs are distinguished by the length spectra of the geodesics in the hyperbolic structure, using the technique defined in~\cite{hodgson1994symmetries}. For some of the manifolds, the 64 bit precision used in the standard version of SnapPy was not sufficient to compute a Dirchlet domain and length spectrum, but we were successful using the old 2.5.1 68K version of SnapPea which uses 80 bit precision arithmetic.
We then checked the results using a new high precision version of SnapPy developed recently by Marc Culler and Nathan Dunfield, which uses quad doubles. Here each length spectrum calculation (up to length 3, with around 50 decimal places of accuracy) took between 25 and 105 seconds on a MacBook with 2.4 GHz Intel Core i7 processor. Note that the running time is substantially better in some cases as compared with the homology of covers technique used before. 

In our experiments, the shortest geodesic is always one of the core circles added in Dehn filling the link. When this is true the manifolds can also be distinguished by looking at the complement of this shortest geodesic --- these 1-cusped hyperbolic 
3-manifolds can be compared by looking at their canonical decompositions (as computed in SnapPy). 

The last column of Table~\ref{Table:results} gives the complex length $length + i(rotation\   angle)$ for a geodesic distinguishing the two manifolds. The length spectrum gives the complex length for each closed geodesic of length $<2.5$. The command $M.length\_spectrum(3)$ with high precision manifolds was used in SnapPy to find the length spectrum up to length $3$.

\begin{table}
\rowcolors{1}{lightgray}{white}
\caption{{\footnotesize Distinguishing manifolds from $HG8QI$ using volume, covers and complex length of geodesics.}}
\label{Table:results}
\scalebox{0.64} {
\begin{tabular}{ l | l | l | p{7cm}|l | l | p{3.7cm} }
  $HG8QI$ & Blink & Volume & Differing Homology of Covers &{Degree}& Time for covers& Differing Complex lengths\\
  \hline
  \hline
 $14^t_{24}$ & $T[71]$ & 24.807369734 & Z/63 + Z/63 &5& 0m29.01s & 2.4336 + i 0.0000\\
              & $T[79]$ &              & - &5& & {\color{red}2.3054 + i 0.0000} \\
  \hline
  \hline
  $15^t_{16}$ & $T[118]$ & 28.375305555 & Z/229773 &5&2m17.75s & 2.4273 + i 0.2081\\
              & $T[119]$ &              & - &5&                & 2.3631 + i 0.1591 \\
  \hline
  \hline
  $15^t_{19}$ & $T[122]$ & 27.670218370 & Z/2 + Z/2 + Z/15486 &6&639m52.46s & 2.4551 + i 0.0409  \\
              & $T[148]$ &              & - &6&  & 2.2902 + i 0.0671  \\
  \hline
  \hline
  $15^t_{22}$ & $T[128]$ & 27.932219883 & Z/2 + Z/162 &6&125m51.74s & 2.3357 + i 0.3227  \\
              & $T[141]$ &              & - &6&  & 2.2443 + i 0.3042 \\    
  \hline
  \hline
  $16^t_{42}$ & $T[305]$ & 32.908565776 &{\color{red}Z/2 + Z/2 + Z/2 + Z/2 + Z/114}&{\color{red}4}&{\color{red}0m4.98s} & 1.5138 - i 0.1141\\
              & $T[337]$ &              & - &{\color{red}4}&  & 1.3703 - i 0.0485 \\    
  \hline
  \hline  $16^t_{56}$ & $T[320]$ & 29.436263597 & Z/2 + Z/34 + Z &{\color{red}3}&0m2.63s & 0.5954 - i 0.0058  \\
              & $T[357]$ & 29.460315997  & - &{\color{red}3}& & 0.6027 - i 0.0289  \\    
  \hline
  \hline
  $16^t_{140}$ & $T[423]$ & 30.587901596 & Z/2 + Z/2 + Z/4 + Z/4 + Z/641184 + Z (*) &{\color{red}8}&2m15.14s & 2.2487 - i 0.2325  \\
               & $T[444]$ &   & -(*)  &{\color{red}8}& & 2.3656 - i 0.2740  \\    
  \hline
  \hline
  $16^t_{141}$ & $T[424]$ & 30.707487021 & Z/2+Z/4285014 &5&1m18.90s & 2.2657 + i 0.0000 \\
               & $T[468]$ &   & - &5& & 2.3354 + i 0.0000 \\ 
  \hline
  \hline
  $16^t_{142}$ & $T[425]$ & 30.729338019 & Z/2 + Z/2 + Z/168 &{\color{red}6}&233m45.92s &-  \\
               & $T[435]$ &   & - &{\color{red}6}& & 2.4992 + i 0.0000  \\ 
  \hline
  \hline
  $16^t_{149}$ & $T[437]$ & 33.464380115  & -(*)  &{\color{red}8}& 87m45.14s & 0.4297 + i 0.1524 \\
               & $T[464]$ & 30.867341891 & Z/3 + Z/6 + Z/1980 + Z (*) &{\color{red}8}& & 0.4047 + i 0.0054  \\ 
  \hline
  \hline
  $16^t_{233}$ & $T[631]$ & 29.624669407  & Z/4 + Z/4 + Z/36 + Z/36 &5& 1m5.63s & -
 \\
               & $T[633]$ &  & - &5& & 2.2820 + i 0.0000 \\ 
  \hline  
\end{tabular}
}
{\footnotesize (*) Covers from representations onto $PSL(2,7)$ as computed using GAP/Sage.}
\end{table}


\section{Conclusion}
We provide a complete recipe for obtaining the $\kappa_r$-invariant using only blinks. The homeomorphism problem for the oriented 3-manifolds induced by blinks up to 16 edges is
completely solved, as a result of the analysis of homology of coverings and length spectra, and of previous $\kappa_r$-classification reported in~\cite{lins2007blink}~\href{http://arxiv.org/abs/0702057}{\underline{arxiv}} and~\cite{lins2013B}~\href{http://arxiv.org/abs/1305.4540}{\underline{arxiv}}. Our experiments suggest that sometimes the length spectra technique is superior to the homology of covers technique for distinguishing a given pair of hyperbolic 3-manifolds. These results provide evidence of the truth of the following conjecture: 

\begin{conjecture}
Let $B$ and $B'$ be monochromatic 3-connected blinks inducing a 3-manifold, then $B' \in \{B,-B^*\}$. 
\end{conjecture}

{\color{red}It means that, given a mono3c blink $B$ inducing $3$-manifold $M$, if another $3$-manifold $M'$ is also induced by $B$ (or $-B^*$) then $M \equiv M'$ up to orientable homeomorphisms. So each equivalence class of blinks defined up to coin calculus has at most two mono3c blinks ($B$ and $-B^*$).} The conjecture also implies that the oriented homeomorphism problem for the manifolds induced by mono3c blinks can be replaced by the isomorphism problem for blinks. The latter problem is solvable by an $O(n^2\log n)$-algorithm, where $n$ is the number of edges. An important question remain open: which 3-manifolds correspond to the class of 3-connected monochromatic blinks?  

\subsection{Acknowledgments}
We are grateful to Marc Culler and Nathan Dunfield for their helpful discussions and enthusiastic cooperation in finding some of the distinguishing covers presented here.

\bibliographystyle{plain}
\bibliography{bibtexIndex.bib}

%

\newpage
{\bf APPENDIX}

\section{Classifying the $HG8QI$-classes}

The figures show the quantum invariant in polar form where the angle is divided by $\pi$. The number of states for the $\kappa_r$ invariant is represented by $\#sts$.

\subsection{Distinguishing the 2 members of the $HG8QI$-class  $14^t_{24}$}

\includegraphics[width=16cm]{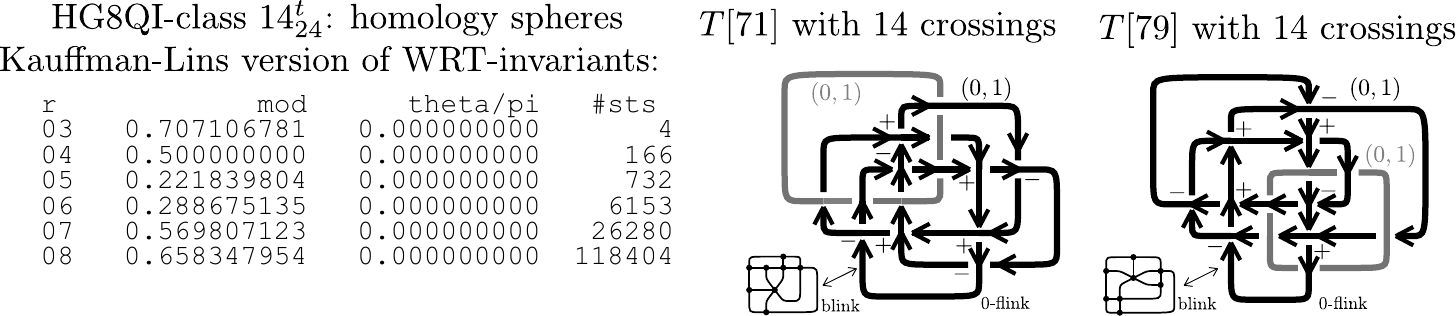} 
Homology of the three 5-covers of $T[71]$:\newline $[Z/3 + Z/3 + Z/3 + Z/3, Z/63 + Z/63, Z/132 + Z/132]$.
\newline\newline
Homology of the three 5-covers of $T[79]$:\newline $[Z/3 + Z/3 + Z/3 + Z/3, Z/213 + Z/213, Z/432 + Z/432]$.
\newline\newline
{Geodesics up to length 2.5 for $T[71]$}:\newline
{\color{red} 1    \qquad 0.4749346632  + 0.0000000000*I}\newline 
1    \qquad 1.9636062749 - 2.7348296878*I \newline
1    \qquad 1.9636062749 + 2.7348296878*I \newline
2    \qquad 2.1697661080 - 1.7487918426*I \newline
2    \qquad 2.1697661080 + 1.7487918426*I \newline
1    \qquad 2.2109848043 - 2.2859284814*I \newline
1    \qquad 2.2109848043 + 2.2859284814*I \newline
{\color{red} 1    \qquad 2.2785556198 + 0.0000000000*I\newline 
1    \qquad 2.4336796074 + 0.0000000000*I *different\newline} 
1    \qquad 2.4886356431 - 2.2055517636*I \newline
\newline
{Geodesics up to length 2.5 for $T[79]$}:\newline
{\color{red} 1    \qquad 0.4749346632 + 0.0000000000*I}\newline 
1    \qquad 1.9636062749 - 2.7348296878*I \newline
1    \qquad 1.9636062749 + 2.7348296878*I \newline
2    \qquad 2.1697661080 - 1.7487918426*I \newline
2    \qquad 2.1697661080 + 1.7487918426*I \newline
1    \qquad 2.2109848043 - 2.2859284814*I \newline
1    \qquad 2.2109848043 + 2.2859284814*I \newline
{\color{red} 1    \qquad 2.2785556198 + 0.0000000000*I\newline 
1    \qquad 2.3053773009 + 0.0000000000*I *different}\newline 
1    \qquad 2.4886356431 - 2.2055517636*I \newline

\subsection{Distinguishing the 2 members of the $HG8QI$-class  $15^t_{16}$}
\includegraphics[width=16.0cm]{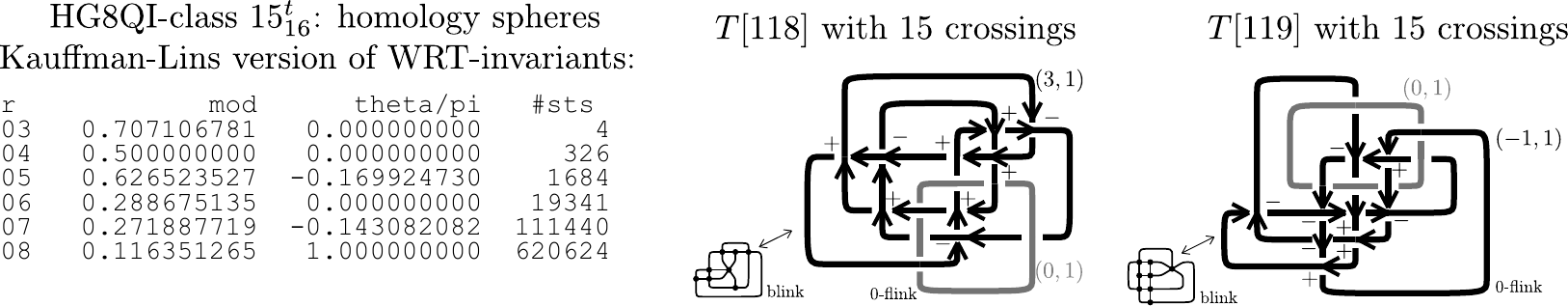} 
Homology of the four 5-covers of $T[118]$:\newline $[Z/229773, Z/1110327, Z/3018207, Z/3699687]$.
\newline\newline
Homology of the four 5-covers of $T[119]$:\newline $[Z/3 + Z/1299909, Z/126627, Z/1052067, Z/4117827]$.
\newline\newline
{Geodesics up to length 2.5 for $T[118]$}:\newline
1    \qquad 0.4253595993 + 0.0518306998*I \newline
1    \qquad 2.0621344966 + 2.6433292590*I \newline
2    \qquad 2.1480316307 - 1.7525028108*I \newline
2    \qquad 2.1718731334 + 1.7708601142*I \newline
1    \qquad 2.2202593862 + 2.7966639551*I \newline
1    \qquad 2.2846442364 + 2.1057757387*I \newline
1    \qquad 2.4273730382 + 0.2081689745*I *different\newline
1    \qquad 2.4489374482 - 1.9955471274*I \newline
{\color{red} 1    \qquad 2.4880484482 + 2.4037153344*I \newline
1    \qquad 2.4897680693 - 1.7075203355*I \newline}
\newline
{Geodesics up to length 2.5 for $T[119]$}:\newline
1    \qquad 0.4253595993 + 0.0518306998*I \newline
1    \qquad 2.0621344966 + 2.6433292590*I \newline
2    \qquad 2.1480316307 - 1.7525028108*I \newline
2    \qquad 2.1718731334 + 1.7708601142*I \newline
1    \qquad 2.2202593862 + 2.7966639551*I \newline
1    \qquad 2.2846442364 + 2.1057757387*I \newline
1    \qquad 2.3631432567 + 0.1591389238*I *different\newline
1    \qquad 2.4489374482 - 1.9955471274*I \newline
{\color{red}  1    \qquad 2.4880484482 + 2.4037153344*I \newline
 1    \qquad 2.4897680693 - 1.7075203355*I \newline}
\newpage
\subsection{Distinguishing the 2 members of the $HG8QI$-class  $15^t_{19}$}
The $HG8QI$-class $15_{19}$ provided us 
with an example of the toughness of the computation
involved to obtain the coverings. It took SnapPy/Sage more than 7 hours
to obtain the (reported) 20 coverings.

\includegraphics[width=12.5cm]{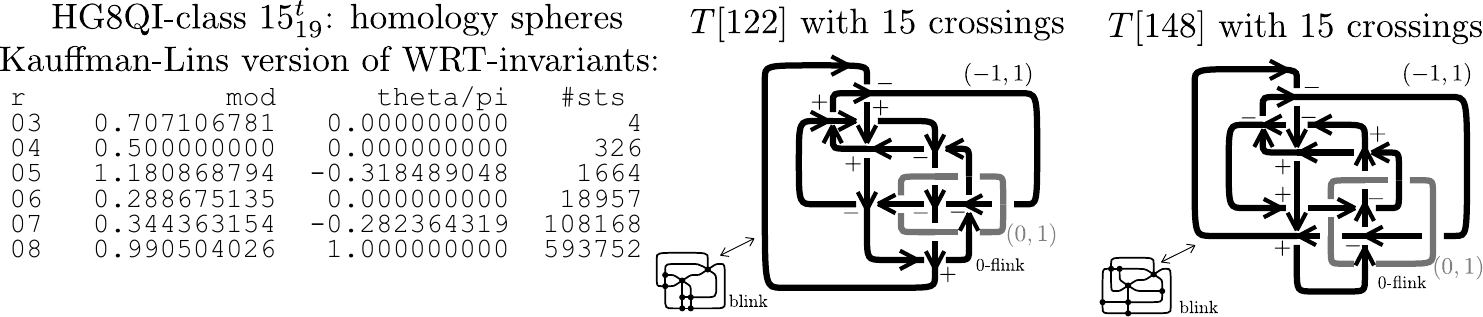} 
\newline
Homology of the twenty $6$-covers of
$T[122]$:\newline
$[
 Z/2 + Z/2 + Z/15486,
 Z/2 + Z/2 + Z/15486,
 Z/2 + Z/894,
 Z/2 + Z/894,
 Z/3 + Z/3 + Z/732,
 Z/3 + Z/3 + Z/732,
 Z/3 + Z/3 + Z/6552,
 Z/3 + Z/3 + Z/6552,
 Z/4 + Z/12,
 Z/4 + Z/12,
 Z/6 + Z/402,
 Z/6 + Z/402,
 Z/24 + Z,
 Z/24 + Z,
 Z/531,
 Z/531,
 Z/549,
 Z/549,
 Z/4716, \\
 Z/4716]$.
 \newline\newline
Homology of the twenty 6-covers of $T[148]$: \newline
$[Z/2 + Z/2 + Z/114,
 Z/2 + Z/2 + Z/384,
 Z/2 + Z/966,
 Z/3 + Z/3 + Z/51,
 Z/3 + Z/3 + Z/1302,
 Z/3 + Z/24 + Z,
 Z/6 + Z/6 + Z/444,
 Z/6 + Z/138,
 Z/6 + Z/150,
 Z/6 + Z/162,
 Z/6 + Z/450,
 Z/6 + Z/4842,
 Z/8 + Z/24,
 Z/18 + Z/450,
 Z/1341,
 Z/4356,
 Z/4866,
 Z/11025,
 Z/66402,\\
 Z/71496]$.
\newline\newline
{Geodesics up to length 2.5 for $T[122]$}:\newline
1    \qquad 0.4366419159 + 0.0557439922*I \newline
1    \qquad 1.8726428246 + 3.1332895508*I \newline
2    \qquad 2.1487104786 - 1.7490483585*I \newline
2    \qquad 2.1750738046 + 1.7692995461*I \newline
1    \qquad 2.1775924323 - 2.4753690206*I \newline
2    \qquad 2.3574941351 + 2.8080309585*I \newline
1    \qquad 2.3740831595 + 2.0176658347*I \newline
2    \qquad 2.3750076185 + 2.2180532874*I \newline
1    \qquad 2.4551713168 + 0.0409451012*I *different\newline
1    \qquad 2.4581260393 - 1.8172608337*I \newline
\newline
{Geodesics up to length 2.5 for $T[148]$}:\newline
1    \qquad 0.4366419159 + 0.0557439922*I \newline
1    \qquad 1.8726428246 + 3.1332895508*I \newline
2    \qquad 2.1487104786 - 1.7490483585*I \newline
2    \qquad 2.1750738046 + 1.7692995461*I \newline
1    \qquad 2.1775924323 - 2.4753690206*I \newline
1    \qquad 2.2902578695 + 0.0671324677*I *different \newline
2    \qquad 2.3574941351 + 2.8080309585*I \newline
1    \qquad 2.3740831595 + 2.0176658347*I \newline
2    \qquad 2.3750076185 + 2.2180532874*I \newline
1    \qquad 2.4581260393 - 1.8172608337*I \newline
\subsection{Distinguishing the 2 members of the $HG8QI$-class  $15^t_{22}$}
\includegraphics[width=13.0cm]{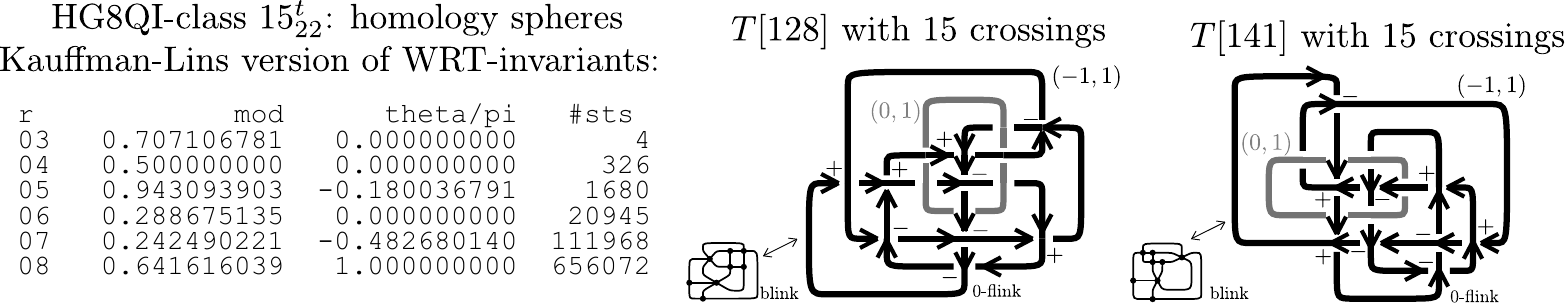} 
\newline
Homology of the fifteen 6-covers of $T[128]$: \newline $[Z/2 + Z/162,
 Z/2 + Z/234,
 Z/2 + Z/414,
 Z/2 + Z/1086,
 Z/3 + Z/3 + Z,
 Z/3 + Z/6 + Z + Z,
 Z/3 + Z/6 + Z + Z,
 Z/3 + Z/6 + Z/18 + Z/720,
 Z/3 + Z/6 + Z/18 + Z/720,
 Z/6 + Z/6 + Z/1296,
 Z/6 + Z/6 + Z/10800,
 Z/9,
 Z/207,
 Z/5571,
 Z/5637]$.
 \newline\newline
Homology of the fifteen 6-covers of $T[141]$: \newline
$[Z/2 + Z/2394,
 Z/3 + Z/3 + Z,
 Z/3 + Z/6 + Z + Z,
 Z/3 + Z/6 + Z + Z,
 Z/3 + Z/6 + Z/18 + Z/720,
 Z/3 + Z/6 + Z/18 + Z/720,
 Z/3 + Z/2439,
 Z/3 + Z/3921,
 Z/6 + Z/6 + Z/144,
 Z/6 + Z/6 + Z/5328,
 Z/6 + Z/1386,
 Z/6 + Z/1482,
 Z/9,
 Z/207,
 Z/2115]
$.
\newline\newline
{Geodesics up to length 2.5 for $T[128]$}:\newline
1    \qquad 0.4334559667 + 0.0446914297*I \newline
1    \qquad 2.0300279130 + 3.0298962814*I \newline
2    \qquad 2.1510507520 - 1.7515439806*I \newline
2    \qquad 2.1720269840 + 1.7676637154*I \newline
1    \qquad 2.2686413408 - 2.0858848833*I \newline
1    \qquad 2.2799180590 + 2.2333530289*I \newline
1    \qquad 2.3152670340 + 2.5898886981*I \newline
1    \qquad 2.3357228621 + 0.3227678337*I *different \newline
1    \qquad 2.4135767557 - 2.5092704258*I \newline
{\color{red} 1    2.4812726701 + 2.2439657392*I }\newline  
\newline
{Geodesics up to length 2.5 for $T[141]$}:\newline
1    \qquad 0.4334559667 + 0.0446914297*I \newline
1    \qquad 2.0300279130 + 3.0298962814*I \newline
2    \qquad 2.1510507520 - 1.7515439806*I \newline
2    \qquad 2.1720269840 + 1.7676637154*I \newline
1    \qquad 2.2443248232 + 0.3042726024*I *different \newline
1    \qquad 2.2686413408 - 2.0858848833*I \newline
1    \qquad 2.2799180590 + 2.2333530289*I \newline
1    \qquad 2.3152670340 + 2.5898886981*I \newline
1    \qquad 2.4135767557 - 2.5092704258*I \newline
{\color{red} 1    2.4812726701 + 2.2439657392*I} \newline 
\subsection{Distinguishing the 2 members of the $HG8QI$-class  $16^t_{42}$}
\includegraphics[width=16.0cm]{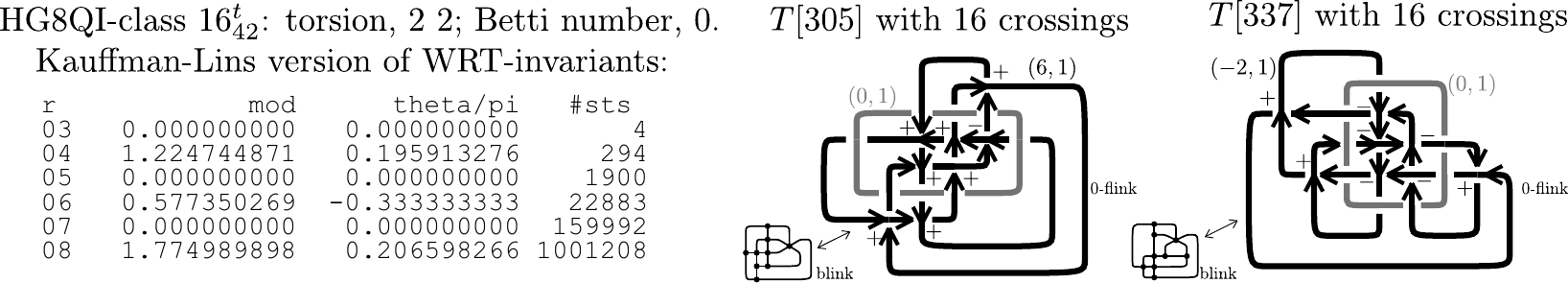} 
\newline
{\color{red}
Homology of the six $4$-covers of $T[305]$: $ \newline
[Z/2 + Z/2 + Z/2 + Z/2 + Z/114,
 Z/2 + Z/29146,
 Z/2 + Z/2 + Z/2 + Z/2 + Z/6,
 Z/2 + Z/2 + Z/2 + Z/2 + Z/6,
 Z/2 + Z/1305434,
 Z/16169 + Z]
  $.
  \newline\newline
Homology of the six $4$-covers of $T[337]$: \newline
$
 [Z/2 + Z/2 + Z/2 + Z/6 + Z/6,
 Z/2 + Z/2 + Z/2 + Z/2 + Z/150,
 Z/2 + Z/2 + Z/2 + Z/2 + Z/150,
 Z/2 + Z/29146,
 Z/2 + Z/1305434,
 Z/16169 + Z]
$.}
\newline\newline
{Geodesics up to length 2.5 for $T[305]$}:\newline
1    \qquad 0.4936558944 + 0.0093839712*I \newline
1    \qquad 1.5138292371 - 0.1141582766*I *different \newline
1    \qquad 2.3192143364 + 2.4343622433*I \newline
{\color{red}  2    \qquad 2.4168583088 + 2.6008225375*I \newline
 2    \qquad 2.4478072620 + 2.0754914062*I \newline
 2    \qquad 2.4969436185 - 2.3591287110*I }\newline 
\newline
{Geodesics up to length 2.5 for $T[337]$}:\newline
1    \qquad 0.4936558944 + 0.0093839712*I \newline
1    \qquad 1.3703361187 - 0.0485212104*I *different\newline
1    \qquad 2.3192143364 + 2.4343622433*I \newline
{\color{red} 2    \qquad 2.4182266342 - 1.8708534578*I \newline
2    \qquad 2.4478072620 + 2.0754914062*I} \newline 
\subsection{Distinguishing the 2 members of the $HG8QI$-class  $16^t_{56}$}
\includegraphics[width=16.0cm]{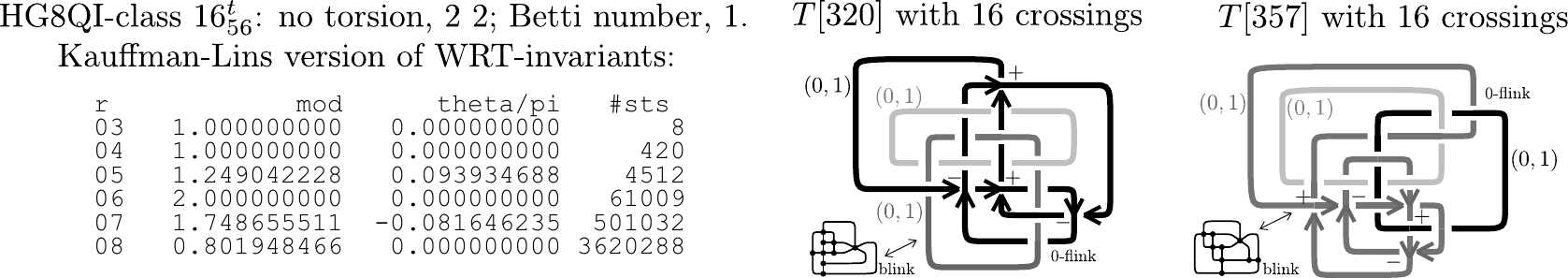} 
\newline
Homology of the five {\color{red}3}-covers of $T[320]$: \newline
$
 [Z, Z/2 + Z/34 + Z, Z/7 + Z/7 + Z, Z/37 + Z, Z/71 + Z]
 $.
 \newline\newline
Homology of the five {\color{red}3}-covers of $T[357]$: \newline
$
[Z, Z/2 + Z/2 + Z/2 + Z/4 + Z, Z/2 + Z/32 + Z, Z/2 + Z/40 + Z, Z/7 + Z/7 + Z]
$.
\newline\newline{Geodesics up to length 2.5 for $T[320]$}:\newline
1    \qquad 0.5954044521 - 0.0058191225*I *different\newline
1    \qquad 1.1282533194 - 0.0471093000*I \newline
{\color{red} 1    \qquad 2.1443810009 + 0.2324419896*I \newline
 1    \qquad 2.2453027917 + 2.4790067843*I \newline
 2    \qquad 2.2945638030 + 1.5857450662*I \newline
 1    \qquad 2.3193876277 - 1.9371730149*I \newline
 2    \qquad 2.3226380083 - 1.6133688144*I \newline
 1    \qquad 2.4011489573 - 2.4082193001*I \newline
 1    \qquad 2.4752251131 - 2.3101512250*I} \newline 
\newline
Geodesics up to length 2.5 for $T[357]$:\newline
%
1    \qquad 0.6027338612 - 0.0289924597*I *different\newline
1    \qquad 1.1061845010 + 0.0327135752*I \newline
{\color{red}  1    2.0359905664 + 0.1785892532*I \newline
 1    \qquad 2.2674125739 - 2.4643790771*I \newline
 1    \qquad 2.3004290474 - 1.9303749709*I \newline
 2    \qquad 2.3005780288 - 1.5973020783*I \newline
 2    \qquad 2.3056141611 + 1.6104156691*I \newline
 1    \qquad 2.3743834017 + 2.4176038653*I \newline
 1    \qquad 2.4385020812 + 2.0093900381*I \newline} 
\subsection{Distinguishing the 2 members of the $HG8QI$-class  $16^t_{140}$}
\includegraphics[width=16.0cm]{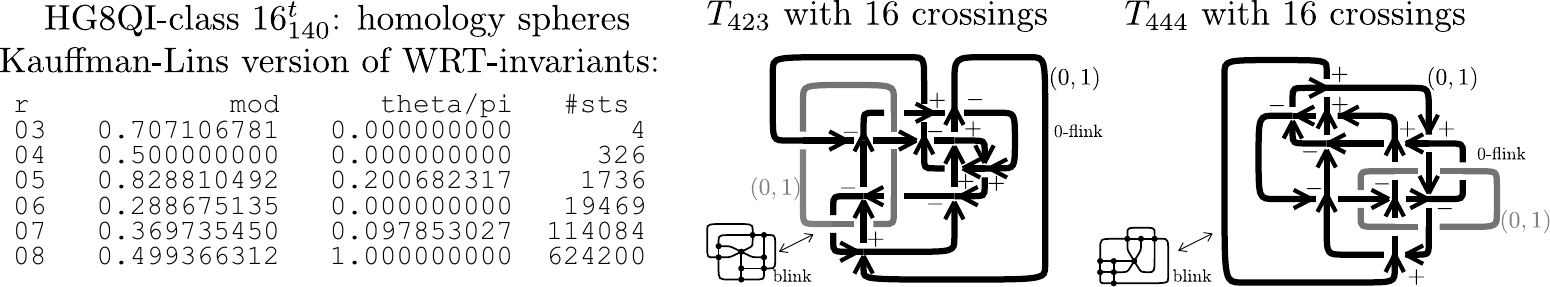} 
\newline
Trying to distinguish these manifolds  
using homology of covers
was initially inconclusive: there are no $k$-covers for $k=2,\ldots, 6$ and the search for
7-covers took more than a week without ending using SnapPy/Sage. 
But we could distinguish using GAP/Sage by examining the special 8-coverings corresponding to homomorphisms onto $PSL(2,{\mathbb F}_7)$, as described in section \ref{subsec:covers} above.
\newline\newline
Homology of  {\color{red}8}-covers of $T[423]$ {\color{red} from homomorphisms onto $PSL(2,{\mathbb F}_7)$}: \newline
$[Z/2 + Z/2 + Z/4 + Z/4 + Z/641184, Z/2 + Z/2 + Z/4 + Z/4 + Z/641184]
$.
\newline\newline
Homology of  {\color{red}8}-covers of $T[444]$ {\color{red} from homomorphisms onto $PSL(2,{\mathbb F}_7)$}: \newline
$[Z/2 + Z/2 + Z/4 + Z/4 + Z/221232, Z/2 + Z/2 + Z/4 + Z/4 + Z/774480]
$.
\newline\newline
{Geodesics up to length 2.5 for $T[423]$}:\newline
1    \qquad 0.4099055382 + 0.0096250743*I \newline
2    \qquad 2.1556583550 - 1.7627803125*I \newline
2    \qquad 2.1599175423 + 1.7660683108*I \newline
1    \qquad 2.2487615848 - 0.2325752662*I *different\newline
2    \qquad 2.2585377651 + 2.8847589106*I \newline
2    \qquad 2.2668792956 - 2.6921566972*I \newline
2    \qquad 2.3701274014 + 2.0417994606*I \newline
{\color{red} 1    2.4208756655 + 2.1803553102*I \newline
1    2.4834089621 + 2.5709250411*I \newline
2    2.4903368368 - 1.8249803819*I }\newline 
 \newline
{Geodesics up to length 2.5 for $T[444]$}:\newline
%
1    \qquad 0.4099055382 + 0.0096250743*I \newline
2    \qquad 2.1556583550 - 1.7627803125*I \newline
2    \qquad 2.1599175423 + 1.7660683108*I \newline
2    \qquad 2.2585377651 + 2.8847589106*I \newline
2    \qquad 2.2668792956 - 2.6921566972*I \newline
1    \qquad 2.3656076568 - 0.2740530565*I *different\newline
2    \qquad 2.3701274014 + 2.0417994606*I \newline
{\color{red} 1    \qquad 2.4208756655 + 2.1803553102*I \newline
1    \qquad 2.4834089621 + 2.5709250411*I \newline
2    \qquad 2.4903368368 - 1.8249803819*I }\newline 
\subsection{Distinguishing the 2 members of the $HG8QI$-class  $16^t_{141}$}
\includegraphics[width=16.0cm]{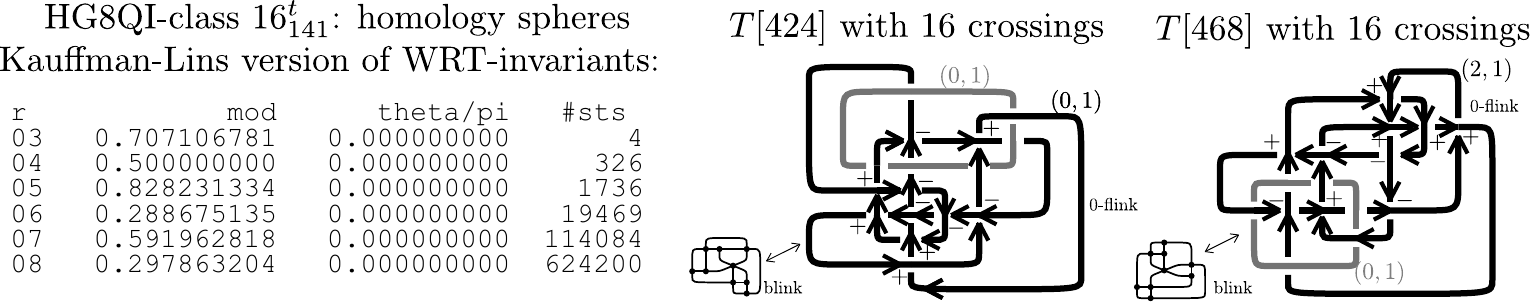} 
\newline
Homology of the three 5-covers of T[424] : \newline
$
[Z/2 + Z/4285014, Z/2 + Z/4285014, Z/3 + Z/3 + Z/3 + Z/3]
$.
\newline\newline
Homology of the three 5-covers of T[468] : \newline
$
[Z/2 + Z/3357564, Z/2 + Z/3357564, Z/3 + Z/3 + Z/3 + Z/3]
$.
\newline
\newline
{Geodesics up to length 2.5 for $T[424]$}:\newline
{\color{red}1    \qquad 0.4070520240 + 0.0000000000*I}\newline 
2    \qquad 2.1573209680 - 1.7650269243*I \newline
2    \qquad 2.1573209680 + 1.7650269243*I \newline
1    \qquad 2.1871322399 - 2.7289199577*I \newline
1    \qquad 2.1871322399 + 2.7289199577*I \newline
{\color{red}1    \qquad 2.2657032372 + 0.0000000000*I  *different }\newline 
1    \qquad 2.3167156129 - 2.0729227053*I \newline
1    \qquad 2.3167156129 + 2.0729227053*I \newline
1    \qquad 2.4360235153 - 2.7467785138*I \newline
{\color{red}1    \qquad 2.4360235153 + 2.7467785138*I \newline
1    \qquad 2.4940537187 - 1.7697279911*I \newline
1    \qquad 2.4940537187 + 1.7697279911*I \newline}
\newline
{Geodesics up to length 2.5 for $T[468]$}:\newline
{\color{red} 1    \qquad 0.4070520240 + 0.0000000000*I}\newline 
2    \qquad 2.1573209680 - 1.7650269243*I \newline
2    \qquad 2.1573209680 + 1.7650269243*I \newline
1    \qquad 2.1871322399 - 2.7289199577*I \newline
1    \qquad 2.1871322399 + 2.7289199577*I \newline
1    \qquad 2.3167156129 - 2.0729227053*I \newline
1    \qquad 2.3167156129 + 2.0729227053*I \newline
{\color{red} 1    \qquad 2.3354490196 + 0.0000000000*I}\newline 
1    \qquad 2.4360235153 - 2.7467785138*I \newline
{\color{red} 1    \qquad 2.4360235153 + 2.7467785138*I \newline
1    \qquad 2.4940537187 - 1.7697279911*I \newline
1    \qquad 2.4940537187 + 1.7697279911*I \newline}
\subsection{Distinguishing the 2 members of the $HG8QI$-class  $16^t_{142}$}
\includegraphics[width=16.0cm]{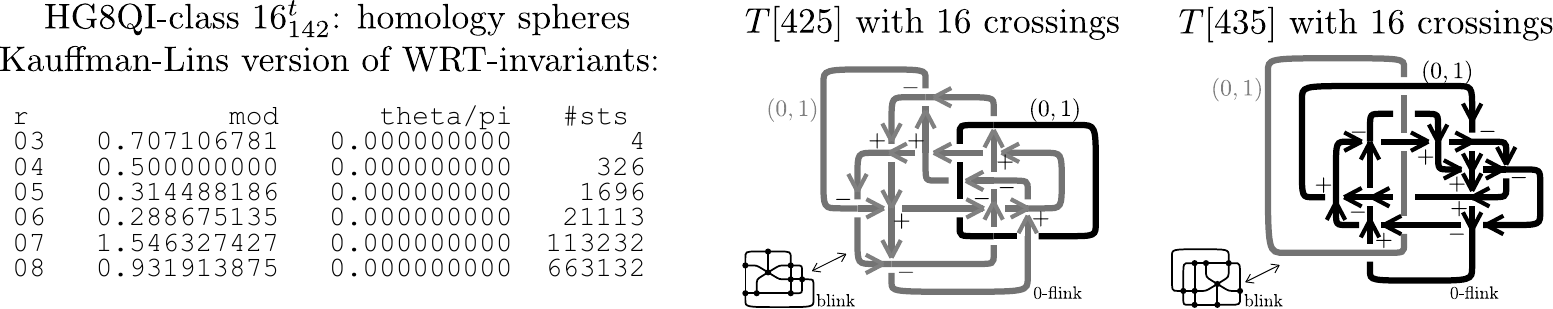} 
\newline
Homology of the twenty nine {\color{red}6}-covers of {\color{red}$T[425]$}: \newline
$
[Z/2 + Z/2 + Z/2 + Z/2 + Z/4,
 Z/2 + Z/2 + Z/2 + Z/2 + Z/4,
 Z/2 + Z/2 + Z/2 + Z/2 + Z/12 + Z/12 + Z,
 Z/2 + Z/2 + Z/168,
 Z/2 + Z/2 + Z/168,
 Z/2 + Z/2 + Z/12804,
 Z/2 + Z/2 + Z/12804,
 Z/2 + Z/2 + Z/42006,
 Z/2 + Z/2 + Z/42006,
 Z/3 + Z/3 + Z,
 Z/3 + Z/3 + Z,
 Z/3 + Z/3 + Z/162,
 Z/3 + Z/3 + Z/162,
 Z/3 + Z/3 + Z/1668,
 Z/3 + Z/3 + Z/1668,
 Z/3 + Z/84,
 Z/3 + Z/84,
 Z/3 + Z/990,
 Z/3 + Z/990,
 Z/3 + Z/13548,
 Z/3 + Z/13548,
 Z/3 + Z/26358,
 Z/3 + Z/26358,
 Z/4 + Z/4 + Z,
 Z/4 + Z/4 + Z,
 Z/6 + Z/6 + Z,
 Z/57 + Z/57 + Z,
 Z/7392,
 Z/7392]
$.
\newline\newline
Homology of the twenty nine {\color{red}6}-covers of $T[435]$: \newline
$
[Z/2 + Z/2 + Z/2 + Z/2 + Z/4,
 Z/2 + Z/2 + Z/2 + Z/2 + Z/4,
 Z/2 + Z/2 + Z/2 + Z/2 + Z/12 + Z/12 + Z,
 Z/3 + Z/3 + Z,
 Z/3 + Z/3 + Z,
 Z/3 + Z/468,
 Z/3 + Z/468,
 Z/3 + Z/1632,
 Z/3 + Z/1632,
 Z/3 + Z/3618,
 Z/3 + Z/3618,
 Z/3 + Z/47694,
 Z/3 + Z/47694,
 Z/4 + Z/4 + Z,
 Z/4 + Z/4 + Z,
 Z/6 + Z/6 + Z/6324,
 Z/6 + Z/6 + Z/6324,
 Z/6 + Z/6 + Z/61254,
 Z/6 + Z/6 + Z/61254,
 Z/9 + Z/9 + Z,
 Z/126 + Z/126 + Z,
 Z/4200,
 Z/4200,
 Z/10422,
 Z/10422,
 Z/13344,\\
 Z/13344,
 Z/14454,
 Z/14454]
$.
\newline\newline
{Geodesics up to length 2.5 for $T[425]$}:\newline
{\color{red} 1    \qquad 0.4085230467 + 0.0000000000*I}\newline 
1    \qquad 1.9630893605 - 2.7482231246*I \newline
1    \qquad 1.9630893605 + 2.7482231246*I \newline
1    \qquad 1.9774488886 - 2.4603249678*I \newline
1    \qquad 1.9774488886 + 2.4603249678*I \newline
1    \qquad 2.2096915615 - 2.3044648228*I \newline
1    \qquad 2.2096915615 + 2.3044648228*I \newline
{\color{red} 1    \qquad 2.2802137771 + 0.0000000000*I}\newline 
2    \qquad 2.4719187872 - 1.9563368724*I \newline
2    \qquad 2.4719187872 + 1.9563368724*I \newline
1    \qquad 2.4843119704 - 2.2237830532*I \newline
1    \qquad 2.4843119704 + 2.2237830532*I \newline
\newline
{Geodesics up to length 2.5 for $T[435]$}:\newline
%
{\color{red} 1    \qquad 0.4085230467 + 0.0000000000*I}\newline 
1    \qquad 1.9630893605 - 2.7482231246*I \newline
1    \qquad 1.9630893605 + 2.7482231246*I \newline
1    \qquad 1.9774488886 - 2.4603249678*I \newline
1    \qquad 1.9774488886 + 2.4603249678*I \newline
1    \qquad 2.2096915615 - 2.3044648228*I \newline
1    \qquad 2.2096915615 + 2.3044648228*I \newline
{\color{red} 1    \qquad 2.2802137771 + 0.0000000000*I}\newline 
2    \qquad 2.4719187872 - 1.9563368724*I \newline
2    \qquad 2.4719187872 + 1.9563368724*I \newline
1    \qquad 2.4843119704 - 2.2237830532*I \newline
1    \qquad 2.4843119704 + 2.2237830532*I \newline
{\color{red} 1    \qquad 2.4992914148 + 0.0000000000*I  *different}\newline 
\subsection{Distinguishing the 2 members of the $HG8QI$-class  $16^t_{149}$}
\includegraphics[width=16.0cm]{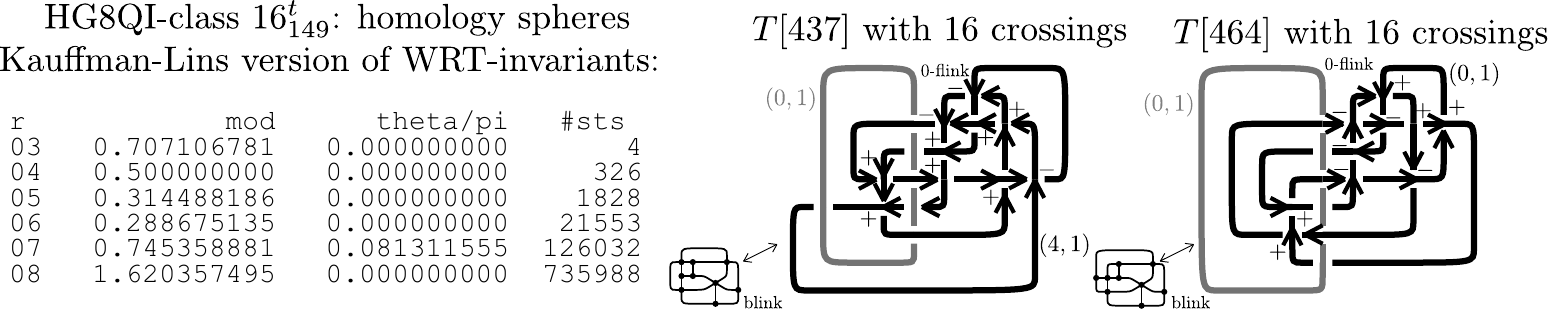} 
We were not able to find covers with SnapPy/Sage after waiting several 
days. However,
the volumes are distinctly different: 33.464380115 for $T[437]$ and 30.8673418910 for $T[464]$.
{\color{red}Using GAP/Sage, we could compute
the special 8-coverings corresponding to homomorphisms onto $PSL(2,{\mathbb F}_7)$}, as described in section~\ref{subsec:covers} above. 
\newline\newline
{\color{red}Homology of the 8-covers} of $T[437]$ {\color{red} from homomorphisms onto $PSL(2,{\mathbb F}_7)$}: \newline
$
[]
$.
\newline\newline
{\color{red}Homology of the 8-covers} of $T[464]$ {\color{red} from homomorphisms onto $PSL(2,{\mathbb F}_7)$}: \newline
$
[Z/3 + Z/6 + Z/1980 + Z,
 Z/3 + Z/6 + Z/1980 + Z,
 Z/6 + Z/1404 + Z,
 Z/6 + Z/3996 + Z,
 Z/6 + Z/1404 + Z,
 Z/6 + Z/3996 + Z,
 Z/3 + Z/6 + Z/612 + Z,
 Z/3 + Z/6 + Z/612 + Z]
$.
\newline\newline
{Geodesics up to length 2.5 for $T[437]$}:\newline
1    \qquad 0.4297247732 + 0.1524342191*I *different\newline
1    \qquad 0.6441175571 - 1.0065104703*I \newline
{\color{red} 2    \qquad 2.2973709255 + 2.3362322468*I \newline
1    \qquad 2.3046129319 - 2.2510636116*I \newline
1    \qquad 2.3573944592 - 0.0440987127*I \newline
2    \qquad 2.4477465768 + 1.8220219514*I \newline
2    \qquad 2.4484860065 - 2.4171117370*I \newline
2    \qquad 2.4876613998 - 2.6700570910*I }\newline
\newline
{Geodesics up to length 2.5 for $T[464]$}:\newline
%
1    \qquad 0.4047727547 + 0.0054277267*I *different\newline
2    \qquad 2.1557443340 - 1.7646186461*I \newline
{\color{red} 2    \qquad 2.1581144241 + 1.7664502140*I \newline
2    \qquad 2.2635115030 + 2.3208406150*I \newline
1    \qquad 2.2853027510 - 2.2278141294*I \newline
1    \qquad 2.3178063815 + 0.0514266001*I \newline
1    \qquad 2.3278290013 - 0.1579902062*I \newline
2    \qquad 2.3704440305 + 2.7711729686*I \newline
1    \qquad 2.3839026067 - 2.3224837801*I \newline
2    \qquad 2.4227459780 - 2.3427847075*I \newline}
\subsection{Distinguishing the 2 members of the $HG8QI$-class  $16^t_{233}$}
\includegraphics[width=16.0cm]{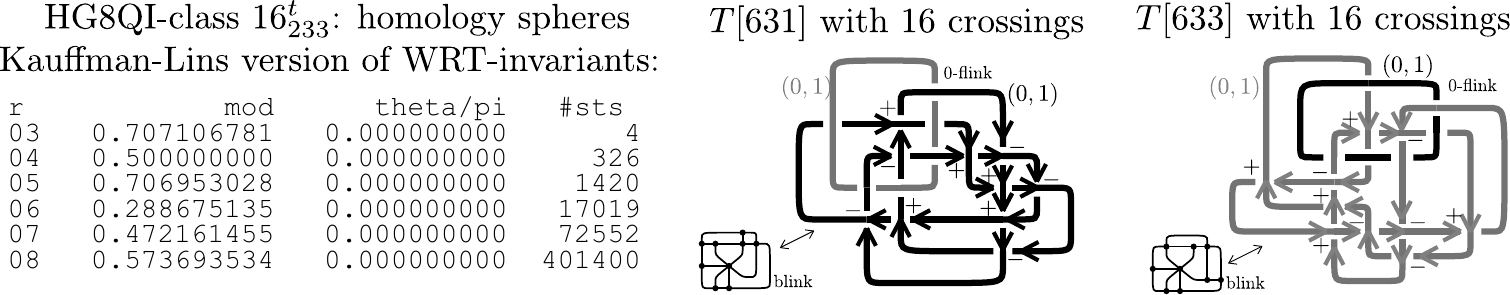} 
\newline
Homology of the two 5-covers of {\color{red}$T[631]$}: \newline
$
[Z/4 + Z/4 + Z/36 + Z/36, Z/21 + Z/21]
$.
\newline\newline
Homology of the two 5-covers of {\color{red}$T[663]$}: \newline
$
[Z/2 + Z/2 + Z/258 + Z/258, Z/171 + Z/171]
$.
\newline\newline
{Geodesics up to length 2.5 for $T[631]$}:\newline
{\color{red}1    \qquad 0.4176144264 + 0.0000000000*I}\newline 
{\color{red}1    \qquad 1.8920374823 + 0.0000000000*I}\newline 
1    \qquad 1.9655039025 - 2.8904656447*I \newline
1    \qquad 1.9655039025 + 2.8904656447*I \newline
2    \qquad 2.1591354297 - 1.7626752539*I \newline
2    \qquad 2.1591354297 + 1.7626752539*I \newline
1    \qquad 2.1603058758 - 2.4909575845*I \newline
1    \qquad 2.1603058758 + 2.4909575845*I \newline
1    \qquad 2.2519402446 - 2.4015681101*I \newline
1    \qquad 2.2519402446 + 2.4015681101*I \newline
1    \qquad 2.4609996482 - 2.3933331803*I \newline
1    \qquad 2.4609996482 + 2.3933331803*I \newline
\newline
{Geodesics up to length 2.5 for $T[663]$}:\newline
%
{\color{red} 1    \qquad 0.4176144264 + 0.0000000000*I}\newline 
{\color{red}1    \qquad 1.8920374823 + 0.0000000000*I}\newline 
1    \qquad 1.9655039025 - 2.8904656447*I \newline
1    \qquad 1.9655039025 + 2.8904656447*I \newline
2    \qquad 2.1591354297 - 1.7626752539*I \newline
2    \qquad 2.1591354297 + 1.7626752539*I \newline
1    \qquad 2.1603058758 - 2.4909575845*I \newline
1    \qquad 2.1603058758 + 2.4909575845*I \newline
1    \qquad 2.2519402446 - 2.4015681101*I \newline
1    \qquad 2.2519402446 + 2.4015681101*I \newline
{\color{red} 1    \qquad 2.2820604976 + 0.0000000000*I *different}\newline
1    \qquad 2.4609996482 - 2.3933331803*I \newline
1    \qquad 2.4609996482 + 2.3933331803*I \newline
\end{document}